\documentclass[12pt]{article}

\usepackage{amsmath}
\usepackage{amsfonts}
\usepackage{amssymb}
\usepackage{graphicx, color}
\usepackage[colorlinks]{hyperref}
\usepackage{epsfig}
\usepackage{epstopdf} 
\usepackage{enumitem}
\usepackage{mathtools}
\usepackage{comment}
\usepackage{soul}

\definecolor{darkgreen}{rgb}{0,0.55,0}

\newtheorem{proposition}{Proposition}[section]
\newtheorem{theorem}{Theorem}[section]
\newtheorem{lemma}[theorem]{Lemma}

\newtheorem{remark}[theorem]{Remark}

\newtheorem{definition}{Definition}

\def\phi{{\varphi}}

\DeclareSymbolFont{AMSb}{U}{msb}{m}{n}
\DeclareMathSymbol{\N}{\mathbin}{AMSb}{"4E}
\DeclareMathSymbol{\Z}{\mathbin}{AMSb}{"5A}
\DeclareMathSymbol{\R}{\mathbin}{AMSb}{"52}
\DeclareMathSymbol{\Q}{\mathbin}{AMSb}{"51}
\DeclareMathSymbol{\I}{\mathbin}{AMSb}{"49}
\DeclareMathSymbol{\C}{\mathbin}{AMSb}{"43}

\DeclareMathOperator*{\esssup}{ess\,sup}

\newcommand{\e}{\varepsilon}

\newcommand{\calA}{{\mathcal A}}
\newcommand{\calH}{{\mathcal H}}

\begin{document}
\title{Existence and structure of solutions for general $P$-area minimizing surfaces}

\author{ {Amir Moradifam \footnote{Department of Mathematics, University of California, Riverside, California, USA. E-mail: amirm@ucr.edu. Amir Moradifam is supported by the NSF grant DMS-1953620.}
}\qquad
Alexander Rowell \footnote{Department of Mathematics, University of California, Riverside, California, USA. E-mail: arowe004@ucr.edu.   } 
}

\date{\today}
\smallbreak \maketitle
\begin{abstract}
We study existence and structure of solutions to the Dirichlet and Neumann boundary problems associated with minimizers of the functional $I(u)=\int_{\Omega} (\phi(x, D u + F)+Hu) \, dx$, where $\phi (x, \xi)$, among other properties,  is convex and homogeneous of degree $1$ with respect to $\xi$. We show that there exists an underlying vector field $N$ that characterizes the existence and structure of  all minimizers.   We also investigate existence of solutions under the barrier condition on $\partial \Omega$.  The results in this paper generalize and unify many results in the literature about existence of minimizers of least gradient problems and $P-$area minimizing surfaces.

\end{abstract}

\section{Introduction and Statement of Results}
In the last two decades, numerous interesting work have been published on existence, uniqueness and regularity of minimizers of functionals of the form 
\[ \int_{\Omega} g(x,Du(x))+k(x,u) \, dx,\] 
where $g$ is convex and $k$ is locally Lipschitz or identically zero.  For background, we encourage the reader to explore the tree of references stemming from \cite{BC,C, CHY, Clar,FT,LM,MR, MT,MT2,MT3, PSTV}.  This paper is a continuation of the authors' work in \cite{MR}, where the authors proved existence and structure of minimizers of P-area minimizing surfaces in the Heisenberg group (see also  \cite{CHY, MR, PSTV} for background literature on P-minmal surfaces in the Heisenberg group).  Let $\Omega$ be a bounded open set in $\R^{2n}$, and 
\[X=(x_1,x'_1,x_2,x'_2, \dots, x_n,x'_n) \in \Omega.\] 
Let  $u: \R^{2n}\rightarrow \R$, and consider the graph $(X,u(X))$ in the Heisenberg group of dimension $2n+1$ with prescribed $p$-mean curvature $H(X)$. Then $u$ satisfies the equation
\begin{equation}\label{meanCurvaturePDE}
\nabla \cdot \left( \frac{\nabla u- X^*}{|\nabla u -X^*|} \right)=H,
\end{equation}
where $X^*=(x'_1,-x_1,x'_2,-x_2, \dots, x'_n, -x_n)$. Equation \eqref{meanCurvaturePDE} is the Euler-Lagrange equation associated  to the energy functional 
\begin{equation}\label{meanCurvatureFunctional}
\mathbb{E}(u)=\int_{\Omega} \left(|\nabla u -X^*|+Hu \right)dx_1 \wedge dx'_1 \wedge \dots \wedge dx_n \wedge dx'_n.
\end{equation}
In \cite{MR} the authors  investigated existence and structure of minimizers of the more general energy functional
\begin{equation}\label{meanCurvatureFunctional2}
\mathbb{I}(u)=\int_{\Omega} \left(a|\nabla u +F|+Hu \right) \, dx,
\end{equation}
under Dirichlet and Neumann boundary conditions and showed that there always exists a vector field $N$ that determines existence and structure of minimizers. Here  $a \in L^{\infty}(\Omega)$ is a positive function and $F \in (L^{\infty}(\Omega))^n$.

In this paper, we study a more general class of functionals which includes \eqref{meanCurvatureFunctional2} as a special case, namely
\begin{equation}\label{mainFunctional}
    I(u)=\int_{\Omega} \phi(x, D u + F)+Hu ,
\end{equation}
where $\phi: \Omega \times \R^{n} \rightarrow \R$ is convex, continuous, and homogeneous function of degree 1 with respect the the second argument.  Unless otherwise stated, we assume that $\Omega$ is a bounded open set in $\R^n$ with Lipschitz boundary, $F \in (L^2(\Omega))^n$, $H\in L^2(\Omega)$, and $\phi$ is assumed to satisfy the following conditions
\begin{enumerate}[label=($C_\arabic*$),leftmargin=2.5\parindent]
    \item  There exists $ \alpha >0$ such that $0 \leq \phi(x, \xi) \leq  \alpha \left| \xi \right|$ for all $\xi \in \R^n$. \label{Condition1}
    \item  $\xi \mapsto \phi(x,  \xi)$ is a norm for every $x$. \label{Condition2} 
\end{enumerate}
While it not generally required, for some of our results we will also assume that 
\begin{enumerate}[label=($C_3$),leftmargin=2.5\parindent]
    \item  There exists $ \beta >0$ such that $0  \leq  \beta \left| \xi \right| \leq \phi(x, \xi)$ for all $\xi \in \R^n$. \label{Condition3}
\end{enumerate}        
This problem is of particular interest since the energy functional $I(u)$ is not strictly convex which makes  analysis of existence and uniqueness of minimizers a highly non-trivial problem. The Rockafellar-Fenchel duality shall play a key role in our study of this problem.  
 
A broad and active area of research is weighted least gradient problems, a special case of (\ref{mainFunctional}) in which  $F\equiv  0$, $H\equiv 0$, and $\phi(x,\xi)=a|\xi|$, where $a\in L^{\infty}(\Omega)$ is a positive function.  This class of sub-class of problems have applications in conductivity imaging and have been extensively studied by many authors, see \cite{HMN, JMN, Mo, Mo1, MNT, MNTa_SIAM, NTT07, NTT08, NTT10, NTT11,  sternberg_ziemer92, sternberg_ziemer, sternbergZiemer93, ST}.  Another interesting special case of (\ref{mainFunctional}) is when $F \equiv 0$, $H\equiv 0$, and $\phi$ is given by 
\[ \phi(x, \xi) = a(x) \left(  \sum_{i,j=1}^{n} \sigma_{0}^{ij} (x) \xi_i \xi_j   \right)^{1/2}, \]

where $\sigma_0=(\sigma^{ij})_{n\times n}$ with $ \sigma^{ij} \in C^{\alpha} (\Omega)$. This problem has applications in imaging of anisotropic conductivity from the knowledge of the interior measurements of current density vector field (see \cite{HMN}).  In \cite{DLPT}, the authors study the case with $H \equiv 0$ and show that, under the so called \textit{bounded slope condition}, the minimizers are Lipschitz continuous. \\

Next we present few preliminaries which are required to understand and the energy functional \eqref{mainFunctional}.  For an arbitrary $u \in BV_{loc}(\R^n)$, an associated measure $\phi(x,Du+F)$ is defined by 
\begin{equation}
    \int_{A} \phi(x, Du+F) = \int_{A} \phi (x, v^u(x)) |Du+F| \hspace{0.5cm} \text{for each bounded Borel set } A,
\end{equation}
with the vector-valued measure $Du+F$ having a corresponding total variation measure $|Du+F|$, and $v^u(x)=\frac{d Du+F}{d |Du+F|}$ is the Radon-Nikodym derivative.  We use standard facts about functions of bounded variation as in \cite{AB}, \cite{JMN}, and \cite{Mo}.  For any open set $U$, we also have 
\begin{equation} \label{generalIBP}
    \int_{U} \phi (x, Du+F)= \sup \left\{  \int_{U} (u \nabla \cdot Y -Y \cdot F) dx : Y\in C_{c}^{\infty}(U;\R^n), \sup \phi^0 (x, Y(x)) \leq 1  \right\},
\end{equation} 
where $\phi(x, \xi)$ has a dual norm on $\R^n$, $\phi^0(x, \xi)$, defined by  
\begin{equation*}
    \phi^0(x, \xi) := \sup \left\{  \xi \cdot p : \phi(x,p) \leq 1  \right\}.
\end{equation*}
As a consequence of condition \ref{Condition1}, the dual norm $\phi^0(x,\cdot)$ has the equivalent definition
\begin{equation}\label{dualNorm}
    \phi^0(x, \xi) = \sup \left\{  \frac{\xi \cdot p}{\phi(x,p) } : p \in \R^n  \right\}.
\end{equation}

\begin{remark}
The definition in (\ref{generalIBP}) allows to define  $\int_{\Omega}\phi(x, Du+F)$ for functions $u\in BV(\Omega)$ with $\nabla u \not \in W^{1,1}(\Omega)$. Indeed the right hand side of (\ref{generalIBP}) is well-defined for any integrable function $u$. To see the motivation behind the definition (\ref{generalIBP}), suppose $u\in W^{1,1}(\Omega)$ and $\phi^0(x,Y) \leq 1$.  For $p=\frac{Du+F}{|Du+F|}$ and $\xi=-Y$ it follows from \eqref{dualNorm} that 
\begin{equation*}
    -Y \cdot \frac{Du+F}{|Du+F|} \leq  \phi\left(x, \frac{Du+F}{|Du+F|} \right).
\end{equation*}
This implies
\begin{align*}
    \int_{\Omega} \phi(x, Du+F) &= \int_{\Omega} \phi\left(x, \frac{Du+F}{|Du+F|} \right) |Du+F| \\
    & \geq \int_{\Omega} -Y \cdot \frac{Du+F}{|Du+F|}|Du+F| \\
    & = \int_{\Omega} -Y \cdot Du - Y \cdot F \\
    & = \int_{\Omega} (u \nabla \cdot Y - Y\cdot F),  \ \   \forall  \ \ Y\in C_{c}^{\infty}(U;\R^n). 
\end{align*}
It is also easy to see that that the inequality above would become an equality in the limit for a sequence of functions $Y_n \in C_{c}^{\infty}(U;\R^n).$
    
\end{remark}

This paper is outlined as follows. In Section 2, we prove existence results under the Neumann boundary condition. In Section 3, we study existence of minimizers with Dirichlet boundary condition.  Finally, in Section 4 we provide existence of P-area minimizing surfaces under a so called \textit{barrier condition} on the boundary $\partial \Omega$. 

\section{Existence of minimizers with Neumann boundary condition} \label{dual}
In this section we study  the minimization problem
\begin{equation}\label{functionalMain0}
	\inf _{u\in \mathring{BV} (\Omega)} I(u):=\int_{\Omega} \phi \left(x, D u + F \right) +Hu,
\end{equation} 
where 
\[\mathring{BV}(\Omega)=\left\{u\in BV(\Omega): \int_{\Omega}u=0 \right\}.\]

We commence our study of minimizers of \eqref{functionalMain0} by applying the Rockefeller-Fenchel duality to the problem. Consider the functions $E: (L^2(\Omega))^n \rightarrow \R $ and $G: \mathring{H}^1(\Omega) \rightarrow \R$ defined as
$$E(b)=\int_{\Omega} \phi \left(x, b+F \right)  \hspace{0.5cm} \text{and} \hspace{0.5cm} G(u)=\int_{\Omega}Hu,$$
where $ \mathring{H}(\Omega)=\{u\in H^1(\Omega): \int_{\Omega}u=0\}$.
Then \eqref{functionalMain0} can be equivalently written as 
\begin{equation}\label{PrimalPronlem0}
	(P) \hspace{0.5cm} \inf_{u \in  \mathring{H}^1(\Omega)} \{  E(\nabla u)+ G(u)\}.  
\end{equation}
The dual problem corresponding to \eqref{PrimalPronlem0}, as defined by Rockafellar-Fenchel duality \cite{ET}, is 
\begin{equation}\label{dualProb0}
(D) \hspace{0.5cm} \max_{b \in (L^{2}(\Omega))^n} \{-E^{*}(b)-G^{*}(-\nabla^{*}b)\}.
\end{equation}
Note that convex functions $E$ and $G$ have convex conjugates $E^*$ and $G^*$.  Furthermore, gradient operator $\nabla: \mathring{H}^1(\Omega) \rightarrow L^2(\Omega)$ has a corresponding adjoint operator $\nabla ^*$.  As computed in \cite{MR}, we have 
\begin{equation*}
	G^*(-\nabla^* b) = \sup_{u \in \mathring{H}^1(\Omega)} \left\{ -\int_{\Omega} \nabla u \cdot b -\int_{\Omega}Hu \right\}.
\end{equation*}  
This can be more explicitly calculated by noting that for all real numbers $c$, $cu \in  \mathring{H}^1 (\Omega)$ whenever $u \in \mathring{H}^1 (\Omega)$.  Thus,
\begin{equation}
G^*(-\nabla^*b)=
	\begin{cases}
    0 & \text{ if } u \in \mathcal{D}_0, \\
	\infty & \text{ if } u \not \in \mathcal{D}_0	
	\end{cases}
\end{equation}
where 
\begin{equation}
    \mathcal{D}_0:=\left \{b\in (L^2(\Omega))^n: \int_{\Omega} \nabla u \cdot b+Hu =0, \ \ \hbox{for all}\ \ u\in  \mathring{H}^1(\Omega)\right \}.
\end{equation}
The computation of $E^*(b)$ is done in Lemma 2.1 of \cite{Mo}, which yields 

\begin{equation}
E^*(b)=
	\begin{cases}
	-\langle F,b \rangle & \text{ if } \phi^0(x, b(x)) \leq  1   \ \ \hbox{in}\ \ \Omega\\
	\infty & \text{ otherwise }.
	\end{cases}
\end{equation}
Thus the dual problem can be rewritten as 
\begin{equation}
   (D) \hspace{0.5cm} \sup \{\langle F,b \rangle: b\in \mathcal{D}_0 \ \ \hbox{and}\ \ \phi^0(x, b(x)) \leq  1   \ \ \hbox{in}\ \ \Omega \}. 
\end{equation}

Let the outer unit normal vector to $\partial\Omega$ be denoted by $\nu_\Omega$. There is a unique function $[b,\nu_\Omega]\in L^{\infty}_{\mathcal{H}^{n-1}}(\partial\Omega)$, whenever $\nabla\cdot b \in L^n(\Omega)$ for every $b \in(L^{\infty}(\Omega))^n$, such that
\begin{equation}\label{trace}
\int_{\partial\Omega}[b,\nu_\Omega]u\,d \mathcal{H}^{n-1}=\int_\Omega u\nabla\cdot b \,dx+\int_\Omega b \cdot D u \, dx,
\quad u\in C^1(\bar\Omega).
\end{equation}
Indeed in \cite{Al, An} it was proved that the integration by parts formula \eqref{trace} holds for every $u\in BV(\Omega)$, as $u\mapsto(b \cdot Du)$ gives rise to a Radon measure on $\Omega$ for $u\in BV(\Omega)$, $b \in(L^{\infty}(\Omega))^n$, and $\nabla \cdot b \in L^n(\Omega)$.

\begin{lemma} \label{LemmaDual}
Let $b \in (L^{\infty}(\Omega))^n \cap \mathcal{D}_0.$ Then 
\[\nabla \cdot b= H-\int_{\Omega}H dx \ \ \hbox{a.e. in}\ \ \Omega,\]
and
\[[b,\nu_{\Omega}]=0 \ \ \mathcal{H}^{n-1}-a.e. \ \ \hbox{ on } \ \ \partial \Omega.\]
\end{lemma}

The above lemma follows directly from equation \eqref{trace} and the definition of $D_0$.  It also provides the insight that every solution $N$ to the dual problem (D) satisfies equation $\nabla \cdot N= H-\int_{\Omega}H dx \ \ \hbox{a.e. in}\ \ \Omega$.  Moreover, at every point on $\partial \Omega$, the unit normal vector is orthogonal to $N$ in a weak sense.  

\begin{theorem} \label{Structure}
Let $\Omega$ be a bounded domain in $\R^n$ with Lipschitz boundary, $F \in (L^2(\Omega))^n$, $H\in L^2(\Omega)$, and $\phi: \Omega \times \R^n \rightarrow \R$ be a convex function satisfying \ref{Condition1} and \ref{Condition2}. Then the duality gap is zero and the dual problem $(D)$ has a solution, i.e. there exists a vector field $N \in \mathcal{D}_0$ with $\phi^0(x,N) \leq 1$ such that 
\begin{equation}\label{dualityGap}
\inf _{u\in \mathring{H}^1(\Omega)} \int_{\Omega} \left( \phi \left(x, D u + F \right) +Hu \right) dx= \langle F, N\rangle.
\end{equation}
Moreover 
\begin{equation}\label{directionParallel}
\phi\left( x, \frac{Du+F}{|Du+F|} \right)= N \cdot \frac{Du+F}{|Du+F|}, \ \ \ \ |Du+F|-a.e. \ \ \hbox{in}\ \ \Omega,
\end{equation}
for any minimizer $u$ of \eqref{PrimalPronlem0}.
\end{theorem}
{\bf Proof.} It is easy to see that  $I(v)=\int_{\Omega} (\phi(x, Dv+F)+Hv)$ is convex, and $J: (L^2(\Omega))^n\rightarrow \R$ with $J(p)=\int_{\Omega} (\phi(x, p+F)+Hu_0) dx$ is continuous at $p=0$, for a fixed $u_0$, due to \ref{Condition2}. Thus, the conditions of Theorem III.4.1 in \cite{ET} are satisfied.  We infer that the optimization problems (D) and (P) have the same optimum value, and the dual problem has a solution $N$ such that the duality gap is zero, i.e.,  \eqref{dualityGap} holds.  

Now let $u\in \mathring{H}^1(\Omega)$ be a minimizer of \eqref{PrimalPronlem0}. Since $N \in \mathcal{D}_0$, we have  
\begin{align*}
    \langle N, F \rangle &= \int_{\Omega} \phi (x, Du+F) +Hu \\
     &= \int_{\Omega} \phi\left( x, \frac{Du+F}{|Du+F|} \right) |Du+F| +\int_{\Omega} Hu \\
    & \geq \int_{\Omega} N \cdot \frac{Du+F}{|Du+F|} |Du+F| +\int_{\Omega} Hu \\
    & = \int_{\Omega} N \cdot (Du+F) + Hu \\
    & = \int_{\Omega}  N \cdot F  +\int_{\Omega} N \cdot Du + Hu. \\
    &= \langle N, F \rangle
\end{align*}
Hence, the inequality above becomes an equality and \eqref{directionParallel} holds. \hfill $\Box$

\begin{remark}
The primal problem $(P)$ may not have a minimizer in $u\in \mathring{H}^1(\Omega)$, but the dual problem $(D)$ always has a solution $N \in (L^2(\Omega))^n$. Note also that the functional $I(u)$ is not strictly convex, and it may have multiple minimizers (see \cite{JMN}). Furthermore, Theorem \ref{Structure} asserts that $N$ determines $\frac{Du+F}{|Du+F|}$, $|Du+F|-$a.e. in $\Omega$, for all minimizers $u$ of $(P)$.  More precisely, since almost everywhere  in $\Omega$ we have 
\[ \phi^0 (x,N) \leq 1 \implies \phi(x, p) \geq N \cdot p  \]
for every $p \in S^{n-1}$.  Therefore, the equality in \eqref{directionParallel} indicates that
\[ \frac{N \cdot p}{\phi (x, p)} \]
is maximized by $p=\frac{Du+F}{|Du+F|}$, $|Du+F|-$a.e.  In the case that $F \equiv 0$, $N$ determines the direction of the gradient of $u$ ($\frac{Du}{|Du|}$) and hence the structure of the level sets of minimizers to $(P)$.
\end{remark}
We proceed to show that a solution  to primal problem $(P)$ exists in $BV(\Omega)$ provided that it is bounded below. The proof that relies on standard facts about $BV$ functions.

\begin{proposition}\label{ExistenceProp} Let $\phi: \Omega \times \R^n \rightarrow \R$ be a convex function satisfying \ref{Condition1}, \ref{Condition2}, and \ref{Condition3}.  If there exists a constant $C$, depending on $\Omega$, such that 

	\begin{equation}\label{BoundExistence}
		\max_{x \in \overline{\Omega}} |H(x)|<C,
	\end{equation}
	then the primal problem (P) has a minimizer.
\end{proposition}
\textbf{Proof.}  Consider the minimizing sequence $u_n$ of functional $I(u)$. By condition \ref{Condition3} we have  
$$ \int_{\Omega} \beta |\nabla u_n + F| + Hu_n \leq \int_{\Omega} \phi(x, \nabla u_n + F) + Hu_n < c,$$
for some constant $c$ independent of $n$.  Moreover, the triangle inquality implies
$$ \int \beta|\nabla u_n | - \int \beta |F| -\int |H| |u_n| \leq \int \beta |\nabla u_n | - \int \beta|F| +\int H u_n \leq \int \beta|\nabla u_n + F| + Hu_n < c$$
and 
$$\int \beta |\nabla u_n | \leq C+ \int |H| |u_n| + \int \beta |F|. $$
Applying Poincar\'e's inequality implies that there exists  a constant $C_{\Omega}$, independent of $n$, where 
$$ \int \beta |\nabla u_n | \leq C+ ||H ||_{L^{\infty}(\Omega)}  C_{\Omega}\int  |\nabla u_n| + \int \beta |F| $$
$$\Rightarrow \left( \beta -C_{\Omega}||H ||_{L^{\infty}(\Omega)}  \right) \int |\nabla u_n| \leq C + \int \beta |F|.  $$
Finally,
$$\int |\nabla u_n | \leq C'= \frac{C+\int \beta |F|}{\left( \beta - C_{\Omega}||H ||_{L^{\infty}(\Omega)}  \right)}$$
provided that $\beta -C_{\Omega}||H ||_{L^{\infty}(\Omega)} >0$ or equivalently 
\[||H ||_{L^{\infty}(\Omega)} \leq C:=\frac{\beta}{C_{\Omega}}.\]

It follows from standard compactness results for $BV$ functions that $u_n$ has a subsequence, denoted by $u_n$ again, such that $u_n $ converges strongly in $L^1$ to a function $\hat{u}\in BV$, and $Du_n$ converges to $D\hat{u}$ in the sense of measures. Since the functional $I(u)$ is lower semicontinuous, $\hat{u}$ is a solution of the primal problem \eqref{functionalMain0}. \hfill $\Box$ 




\section{Existence of minimizers with Dirichlet boundary condition}
Now, let us consider minimizers of the main functional with a given Dirichlet boundary condition on $\partial \Omega$. Let $\Omega$ be a bounded region in $\R^n$ with Lipschitz boundary, $F \in (L^2(\Omega))^n$, $H\in L^2(\Omega)$, $f\in L^1(\partial \Omega)$, $\phi: \Omega \times \R^n \rightarrow \R$ a convex function satisfying \ref{Condition1} and \ref{Condition2}, and the minimization problem becomes
\begin{equation}\label{functionalMain}
	\inf _{u\in BV_f(\Omega)} I(u):=\int_{\Omega} \phi \left(x, D u + F \right) +Hu,
\end{equation} 
where 
\[BV_f(\Omega)=\{u\in BV(\Omega): u|_{\partial \Omega}=f\}.\]
We perform the substitution $\tilde{F}=F+\nabla f$ to rewrite \eqref{functionalMain} in terms of BV functions that are zero on $\partial \Omega$. Since there always exists a function $f \in W^{1,1}(\Omega)$ that is an extension of any function in $L^1(\Omega)$, we have.
\begin{equation*}
	\inf _{u\in BV_0(\Omega)} I(u):=\int_{\Omega} \phi \left(x, D u + \tilde{F} \right) +Hu+\int_{\Omega}Hf dx.
\end{equation*} 
Note that $\int_{\Omega}Hf dx$ is a constant, which implies \eqref{functionalMain} can be represented by the minimization problem
\begin{equation}\label{functionalMainZero}
	\inf _{u\in BV_0(\Omega)} I(u):=\int_{\Omega} \phi \left(x, D u + F \right) +Hu. 
\end{equation} 
In Section \ref{dual} boundedness from below of functional $I(u)$ was sufficient to provide existence of minimizers in $\mathring{BV}(\Omega)$.  This is not the case for \eqref{functionalMain}, nor \eqref{functionalMainZero}.  The main reason for nonexistence of minimizers is that for a given minimizing sequence such that $u_n \rightarrow \hat{u}$ in $L^1(\Omega)$ and $\hat{u} \in BV(\Omega)$, we have
\[I(\hat{u})\leq \inf_{u \in BV_0(\Omega)} I(u),\]
by the lower semicontinuity of $I(u).$  However, since $\partial \Omega$ is a set of measure zero, the trace of $\hat{u}$ is not guaranteed to be zero.  One of our main goals  in this section is to prove existence of minimizers for the highly nontrivial problem \eqref{functionalMainZero}, and in turn for \eqref{functionalMain}.

\subsection{The Dual Problem}

The setup of the dual problem here is identical to that of Section 2, with the exception of the function space of potential solutions.  Let $E: (L^2(\Omega))^n \rightarrow \R $ and $G:H^1_0(\Omega) \rightarrow \R$ be defined as 
$$E(b)=\int_{\Omega} \phi \left(x, b+F \right)  \hspace{0.5cm} \text{and} \hspace{0.5cm} G(u)=\int_{\Omega}Hu.$$
Then \eqref{relaxedProblem} can be equivalently written as 
\begin{equation}\label{PrimalNew}
	(P') \hspace{0.5cm} \inf_{u \in H^1_0(\Omega) } \{  E(\nabla u)+ G(u)\}.  
\end{equation}
The dual problem corresponding to \eqref{PrimalNew}, as defined by Rockafellar-Fenchel duality \cite{ET}, is 
\begin{equation}\label{dualNew}
(D') \hspace{0.5cm} \sup_{b \in (L^{2}(\Omega))^n} \{-E^{*}(b)-G^{*}(-\nabla^{*}b)\}.
\end{equation}
Then $G^*(-\nabla^* b)$ is given by 
\begin{equation*}
	G^*(-\nabla^* b) =\sup_{u \in H^1_0(\Omega) } \left\{ -\int_{\Omega} \nabla u \cdot b -\int_{\Omega}Hu \right\},
\end{equation*}  
and more explicitly 
\begin{equation}
G^*(-\nabla^*b)=
	\begin{cases}
    0 & \text{ if } u \in \widetilde{\mathcal{D}}_0 \\
	\infty & \text{ if } u \not \in \widetilde{\mathcal{D}}_0(\Omega),
	\end{cases}
\end{equation}
where 
\begin{equation}
    \widetilde{\mathcal{D}}_0:=\left \{b\in (L^2(\Omega))^n: \int_{\Omega} \nabla u \cdot b+Hu =0, \ \ \hbox{for all}\ \ u\in H^1_0(\Omega)   \right \}\subseteq \mathcal{D}_0.
\end{equation}
Finally, we use Lemma 2.1 in \cite{Mo} to get
\begin{equation}
E^*(b)=
	\begin{cases}
	-\langle F,b \rangle & \text{ if } \phi^0(x, b(x)) \leq  1   \ \ \hbox{in}\ \ \Omega\\
	\infty & \text{ otherwise }.
	\end{cases}
\end{equation}
We can therefore rewrite the dual problem as
\begin{equation}
   (D') \hspace{0.5cm} \sup \{\langle F,b \rangle: b\in \widetilde{\mathcal{D}}_0 \ \ \hbox{and}\ \ \phi^0(x, b(x)) \leq  1   \ \ \hbox{in}\ \ \Omega \}. 
\end{equation}

A direct application of the integration by parts formula \eqref{trace} implies that $b \in (L^{\infty}(\Omega))^n \cap \widetilde{\mathcal{D}}_0$ if and only if 
\[\nabla \cdot b= H \ \ \hbox{a.e. in}\ \ \Omega.\]\\

Next we proceed to prove the analog of Theorem \ref{Structure}. \\

\begin{theorem} \label{Structure2}
Let $\Omega$ be a bounded domain in $\R^n$ with Lipschitz boundary, $F \in (L^2(\Omega))^n$, $H\in L^2(\Omega)$, $\phi: \Omega \times \R^n \rightarrow \R$ a convex function satisfying \ref{Condition1}, \ref{Condition2}, and assume $(P')$ is bounded below. Then the duality gap is zero and the dual problem $(D')$ has a solution, i.e. there exists a vector field $N \in \widetilde{\mathcal{D}}_0$ with $\phi^0(x,N) \leq 1$ such that 
\begin{equation}\label{dualityGap2}
\inf _{u\in H_0^1(\Omega)} \int_{\Omega} \left( \phi \left(x, D u + F \right) +Hu \right) dx= \langle F, N\rangle.
\end{equation}
Moreover 
\begin{equation}\label{directionParallel2}
\phi\left( x, \frac{Du+F}{|Du+F|} \right)= N \cdot \frac{Du+F}{|Du+F|}, \ \ \ \ |Du+F|-a.e. \ \ \hbox{in}\ \ \Omega,
\end{equation}
for any minimizer $u$ of \eqref{PrimalNew}.
\end{theorem}
{\bf Proof.}  It is easy to show  that $I(v)=\int_{\Omega} (\phi(x, Dv+F)+Hv)$ is convex, and $J: (L^2(\Omega))^n\rightarrow \R$ with $J(p)=\int_{\Omega} (\phi(x, p+F)+Hu_0) dx$ is continuous at $p=0$, for a fixed $u_0$, due to \ref{Condition2}. Thus, the conditions of Theorem III.4.1 in \cite{ET} are satisfied.  We infer that the optimal values of $(D)$ and (P) are equal, and the dual problem has a solution $N$ such that the duality gap is zero, i.e. \eqref{dualityGap2} holds.  

Now let $u\in A_0$ be a minimizer of \eqref{PrimalNew}. Since $N \in \widetilde{\mathcal{D}}_0$ 
\begin{align*}
    \langle N, F \rangle &= \int_{\Omega} \phi (x, Du+F) +Hu \\
     &= \int_{\Omega} \phi\left( x, \frac{Du+F}{|Du+F|} \right) |Du+F| +\int_{\Omega} Hu \\
    & \geq \int_{\Omega} N \cdot \frac{Du+F}{|Du+F|} |Du+F| +\int_{\Omega} Hu \\
    & = \int_{\Omega} N \cdot (Du+F) + Hu \\
    & = \int_{\Omega}  N \cdot F  +\int_{\Omega} N \cdot Du + Hu. \\
    &= \langle N, F \rangle
\end{align*}
Hence the inequality becomes an equality and \eqref{directionParallel2} holds. \hfill $\Box$
\begin{remark}
Similar to the comments we made in Remark 3.1., the primal problem $(P')$ may not have a minimizer in $H_0^1 (\Omega)$, but the dual problem $(D')$ always has a solution $N \in (L^2(\Omega))^n$. Note also that the functional $I(u)$ is not strictly convex, and it may have multiple minimizers (see \cite{JMN}). Furthermore, Theorem \ref{Structure2} asserts that $N$ determines $\frac{Du+F}{|Du+F|}$, $|Du+F|-$a.e. in $\Omega$, for all minimizers $u$ of $(P')$.  See Remark 3.1 for more details. 
\end{remark}

\subsection{The relaxed problem } \label{relaxedSection}
Now we investigate the existence of minimizier for the relaxed problem associated to \eqref{functionalMainZero}, namely
\begin{equation}\label{relaxedProblem}
\inf_{u\in A_0}I(u)= \inf_{u\in A_0} \int_{\Omega}(\phi(x, Du+F)+Hu)dx  +\int_{\partial \Omega} \phi(x, \nu_{\Omega} )|u| ds,
\end{equation}
where 
$$A_0 := \left \{  u \in H^1 (\R^n) : u=0 \text{ in } \Omega^c   \right \}. $$
The benefit of considering the relaxed problem above is that any minimizing sequence of \eqref{relaxedProblem} converges to a minimizer in $A_0$.  This convergence result is not guaranteed for \eqref{functionalMainZero}.  It can be easily verified that Proposition \ref{ExistenceProp} can be adapted to the relaxed problem, and \eqref{functionalMainZero} has a solution in $A_0$ when bounded below.

\begin{proposition} \label{PropLast} Let $\phi: \Omega \times \R^n \rightarrow \R$ be a convex function satisfying \ref{Condition1}, \ref{Condition2}, and \ref{Condition3}.  If there exists a constant $C$, depending on $\Omega$, such that

	\begin{equation}
		\max_{x \in \overline{\Omega}} |H(x)|<C,
	\end{equation}
	then the primal problem \eqref{functionalMainZero} has a minimizer in $A_0$.
\end{proposition}
{\bf Proof.}  Note that $\hat{u}\in A_0$ whenever $u_n \in A_0$ converges to $\hat{u}$ in $L^1(\Omega)$.  Then the proof follows as outlined in Proposition \ref{ExistenceProp}. \hfill $\Box$
\\

The stage is now set for the major result of this section.  While proving existence of minimizers to \eqref{functionalMainZero} is very difficult, the following theorem demonstrates how problems \eqref{functionalMainZero} and \eqref{relaxedProblem} are related.

\begin{theorem}\label{twoPrimalProblems}
	Let $\Omega \subset \R ^n$ be a bounded open set with Lipschitz boundary, $F \in (L^2(\Omega))^n$, $H\in L^2(\Omega)$, and $\phi: \Omega \times \R^n \rightarrow \R$ a convex function satisfying \ref{Condition1}, \ref{Condition2}, and \ref{Condition3}.  If the minimization problem \eqref{functionalMainZero} is bounded below, then

	\begin{equation}\label{MinimizerA0}
		\min_{u \in A_0} \left( \int_{\Omega} (\phi(x, D u +F) +Hu)dx +\int_{\partial \Omega} \phi(x, \nu_{\Omega} )|u| ds  \right) = \inf_{\substack{u\in BV_0(\Omega) }} \int_{\Omega} \phi(x, D u +F) +Hu 
	\end{equation}
Moreover, if $u$ is a minimizer of \eqref{relaxedProblem}, then 
\begin{equation}\label{BoundaryResult1}
\phi(x, \nu_\Omega) =[N, sign(-u)\nu_{\Omega}] \ \ \mathcal{H}^{n-1}-a.e. \ \ \hbox{on}\ \ \partial \Omega. 
\end{equation}

\end{theorem}
{\bf Proof.}  It can be easily shown that $BV_0(\Omega)$ has a continuous embedding into $A_0$, which implies  
$$\min_{u \in A_0} \left( \int_{\Omega} (\phi(x, D u +F) +Hu)dx +\int_{\partial \Omega} \phi(x, \nu_{\Omega} )|u| ds \right) \leq \inf_{\substack{u\in BV_0(\Omega) }} \int_{\Omega} \phi(x, D u +F) +Hu.$$
It follows from Theorem \ref{Structure2} that there exists a vector field $N \in \widetilde{\mathcal{D}}_0$ with 
\[\phi\left( x, \frac{Du+F}{|Du+F|} \right)= N \cdot \frac{Du+F}{|Du+F|}, \ \ \ \ |Du+F|-a.e. \ \ \hbox{in}\ \ \Omega.\]
Consider minimizer $u$ of the relaxed problem with $u|_{\partial \Omega}=h|_{\partial \Omega}$, where $h\in W^{1,1}(\Omega)$. Since $u-h \in \widetilde{\mathcal{D}}_0$, we have  
\begin{align*}
	 \min_{u \in A_0} ( \int_{\Omega} (\phi(x, D u +F) +Hu)dx &+\int_{\partial \Omega} \phi(x, \nu_{\Omega} )|u| ds ) =
    \int_{\Omega} \phi(x, D u +F) +Hu +\int_{\partial \Omega} \phi(x, \nu_{\Omega} )|u|  \\
    &= \int_{\Omega} \phi\left( x, \frac{Du+F}{|Du+F|} \right) |Du+F| +\int_{\Omega} Hu +\int_{\partial \Omega} \phi(x, \nu_{\Omega} )|u|\\
    & \geq \int_{\Omega} N \cdot \frac{Du+F}{|Du+F|} |Du+F| +\int_{\Omega} Hu +\int_{\partial \Omega} \phi(x, \nu_{\Omega} )|u|\\
    & = \int_{\Omega} N \cdot (D u +F)  +Hu +\int_{\partial \Omega} \phi(x, \nu_{\Omega} )|u|\\
    & = \int_{\Omega} N \cdot F + \int_{\Omega} N \cdot D u +Hu +\int_{\partial \Omega} \phi(x, \nu_{\Omega} )|u|\\
    & = \langle N, F \rangle + \int_{\Omega} N \cdot D (u-h) + H(u-h) \\
    & + \int_{\Omega} N \cdot D h +Hh +\int_{\partial \Omega} \phi(x, \nu_{\Omega} )|h|\\
    &= \langle N, F \rangle + \int_{\Omega} N \cdot D h +Hh +\int_{\partial \Omega} \phi(x, \nu_{\Omega} )|h|\\     
    &= \langle N, F \rangle + \int_{\partial \Omega}  [N , \nu_{\Omega}]h  +\int_{\partial \Omega} \phi(x, \nu_{\Omega} )|h|\\     
    & \geq \langle N, F \rangle\\
    & = \inf_{\substack{u\in BV_0(\Omega) }} \int_{\Omega} \phi(x, D u +F) +Hu. 
\end{align*}
The last inequality was achieved using integration by parts and the fact that $\phi^0 (x,N) \leq 1 \implies [N, \nu_{\Omega} ] \leq \phi (x, \nu_{\Omega})$.  Therefore, \eqref{MinimizerA0} holds and all the inequalities in the above computation are equalities. This provides the relationship $\int_{\partial \Omega}  [N , \nu_{\Omega}]h  +\int_{\partial \Omega} \phi(x, \nu_{\Omega} )|h|=0$, which implies that \eqref{BoundaryResult1} holds. \hfill $\Box$
\\

The next theorem follows directly from Theorem \ref{Structure2} and Theorem \ref{twoPrimalProblems}.  
\begin{theorem}\label{summaryTheorem}
Let $\Omega \subset \R ^n$ be a bounded open set with Lipschitz boundary, $F \in (L^2(\Omega))^n$, $H\in L^2(\Omega)$, $\phi: \Omega \times \R^n \rightarrow \R$ a convex function satisfying \ref{Condition1}, \ref{Condition2}, \ref{Condition3}, and assume $(P')$ is bounded below. Then there exists a vector field $N \in \widetilde{\mathcal{D}}_0$ with $\phi^0(x,N) \leq 1$  such that 
\begin{equation}
\phi\left( x, \frac{Du+F}{|Du+F|} \right)= N \cdot \frac{Du+F}{|Du+F|}, \ \ \ \ |Du+F|-a.e. \ \ \hbox{in}\ \ \Omega,
\end{equation}
for any minimizer $u$ of \eqref{functionalMainZero}. Moreover, every minimizer of \eqref{functionalMainZero} is a minimizer of \eqref{relaxedProblem}, and  if $u$ is a minimizer of \eqref{relaxedProblem}, then 
\begin{equation}\label{BoundaryResult}
\phi(x, \nu_\Omega) =[N, sign(-u)\nu_{\Omega}] \ \ \ \mathcal{H}^{n-1}-a.e. \ \ \hbox{on}\ \ \partial \Omega. 
\end{equation}

\end{theorem}

\section{Existence of minimizers under the Barrier condition}
Consider $F\in (L^1(\Omega)^n)$, $H \in L^{\infty} (\Omega)$, and $\psi:\R^n \times BV_0(\Omega)$ given to be 
\begin{equation}\label{AlteredIntegrand}
    \psi(x, u) := \phi(x, Du + F\chi_{E_u}) + Hu,
\end{equation}
with $E_u$ representing the closure of the support of $u$ in $\Omega$.  We also define the $\psi$-perimeter of $E$ in $A$ by
\begin{equation*}
    P_\psi(E;A):= \int_{A} \phi \left(x, D \chi_E + F \chi_E \right) + H \chi_E. 
\end{equation*}

\begin{definition} \label{AreaMinimizingDef} A function $u \in BV(\R^n)$ is $\psi$-total variation minimizing in $\Omega \subset \R^n$ if
    $$
   \int_{\Omega}\psi(x,u) \leq \int_{\Omega}\psi(x,v) \text{  for all } v \in BV(\R^n) \text{ such that } u= v \text{ a.e. in } \Omega^c.
    $$
    
Also a set $E \subset \R ^n$ of finite perimeter is $\psi$-area minimizing in $\Omega$ if
\[P_\psi (E;\Omega) \leq P_\psi (\tilde{E}) \]

  for all $ \tilde{E} \subset \R^n \text{ such that } \tilde{E} \cap \Omega^c = E \cap \Omega ^c \text{ a.e.}$.

\end{definition}

In order to state the two major results of this section, Theorems \ref{SuperLevelThm} and \ref{MainThm}, we need the following preliminary lemmas.  The argument is this section are inspired by and similar to those in \cite{MR}. For a given function $u \in BV(\Omega)$, define functions 
\begin{equation}\label{u1u2}
    u_1= \max (u- \lambda, 0) \text{ and } u_2=u-u_1,
\end{equation}
for an arbitrary $\lambda \in \R$.  Moving forward we shall use the function
\begin{equation}\label{chiEpsilon}
    \chi_{\epsilon, \lambda} := \min \left( 1, \frac{1}{\epsilon} u_1  \right) =
    \begin{cases}
        0 & \text{ if } u \leq \lambda, \\
        \frac{1}{\epsilon} (u-\lambda) & \text{ if } \lambda < u \leq \lambda +\epsilon, \\
        1 & \text{ if } u > \lambda + \epsilon.
    \end{cases}
\end{equation}
which is shown to be $\psi$-total variation minimizing in Theorem \ref{SuperLevelThm}.

\begin{lemma} \label{LSC-Lemma} For $\chi_{\epsilon,\lambda}$ as defined in \eqref{chiEpsilon},
\[ P_{\psi}(E,\Omega) \leq \liminf_{\epsilon \rightarrow 0} \int_{\Omega} \phi(x, D \chi_{\epsilon,\lambda}+F\chi_{\epsilon, \lambda})+H \chi_{\epsilon, \lambda}. \]
\end{lemma}
{\bf Proof.} Due to condition \ref{Condition2} we have 
\begin{eqnarray*}
&& \int_{\Omega} \phi(x, D \chi_{\epsilon,\lambda}+F\chi_{\epsilon, \lambda})+H \chi_{\epsilon, \lambda} -\int_{\Omega} \phi(x, D \chi_{E}+F\chi_{E})+H \chi_{E} \\
&=&  \int_{\Omega \cap \{\lambda-\epsilon<u < \lambda+\epsilon\}} \phi(x, D \chi_{\epsilon,\lambda}+F\chi_{\epsilon, \lambda})+H \chi_{\epsilon, \lambda} - \phi(x, D \chi_{E}+F\chi_{E})-H \chi_{E} \\
&\geq &  \int_{\Omega \cap \{\lambda-\epsilon<u < \lambda+\epsilon\}} \phi(x, D \chi_{\epsilon,\lambda})-\phi(x, F\chi_{\epsilon, \lambda})+H \chi_{\epsilon, \lambda} - \phi(x, D \chi_{E})-\phi(x, F\chi_{E})-H \chi_{E} \\
&=&  \int_{\Omega \cap \{\lambda-\epsilon<u < \lambda+\epsilon\}} \phi(x, D \chi_{\epsilon,\lambda})- \phi(x, D \chi_{E})+H \chi_{\epsilon, \lambda} -H \chi_{E} -\phi(x, F\chi_{\epsilon, \lambda})- \phi(x, F\chi_{E})\\
&=& \int_{\Omega } \phi(x, D \chi_{\epsilon,\lambda})- \int_{\Omega}\phi(x, D \chi_{E})+ \int_{\Omega} (H \chi_{\epsilon, \lambda} -H \chi_{E})\\
&& - \int_{\Omega \cap \{\lambda-\epsilon<u < \lambda+\epsilon\}} \phi(x, F\chi_{\epsilon, \lambda})+ \phi(x, F\chi_{E}).\\
\end{eqnarray*}
Since the last two integrals converge to zero as $\epsilon \rightarrow 0$, 

\begin{eqnarray*}
&&  \liminf_{\epsilon \rightarrow 0} \int_{\Omega} \phi(x, D \chi_{\epsilon,\lambda}+F\chi_{\epsilon, \lambda})+H \chi_{\epsilon, \lambda}- P_{\psi}(E,\Omega)\\
&= & \liminf_{\epsilon \rightarrow 0}\int_{\Omega} \phi(x, D \chi_{\epsilon,\lambda}+F\chi_{\epsilon, \lambda})+H \chi_{\epsilon, \lambda} -\int_{\Omega} \phi(x, D \chi_{E}+F\chi_{E})+H \chi_{E} \\
&\geq & \liminf_{\epsilon \rightarrow 0} \int_{\Omega } \phi(x, D \chi_{\epsilon,\lambda})- \int_{\Omega}\phi(x, D \chi_{E})\geq 0,\\
\end{eqnarray*}
where the lower semi-continuity of $\int_{\Omega}\phi(x, Dv)$ justifies the last inequality (see \cite{JMN}).  \hfill $\Box$

\vspace{.2cm}

The outer and inner trace of $w$ on $\partial \Omega$ are denoted by $w^+$ and $w^{-}$ respectively, under the assumptions that $\Omega$ is an open set with Lipschitz boundary and $w\in BV(\R^n)$.

\begin{lemma}\label{Relaxed-Lemma} Suppose $\Omega\subset \R^n$ is a bounded open region with Lipschitz boundary, $g\in L^1(\partial \Omega ; \mathcal H^{n-1})$, and define
\[
I_\psi(v ; \Omega, g) := \ \int_{\partial \Omega} \phi(x, g - v^- +F_{\chi_v})
d\mathcal H^{n-1}
+
\int_{\Omega}\psi(x,Dv).
\]
Then $u\in BV(\R^n)$ is $\psi$-total variation minimizing in $\Omega$ 
if and only if $u|_\Omega$ minimizes
$I_\psi ( \, \cdot \, ; \Omega, g)$ for some $g$, and moreover $g = u^+$.
\label{lem:reformulate}\end{lemma}

{\bf Proof:}
Note that $v^+, v^- \in L^1(\partial \Omega;\mathcal H^{n-1})$ whenever $v\in BV(\R^n)$.  Conversely, there is a $v\in BV(\R^n)$ with $g = v^+$ for each $g\in L^1(\partial \Omega;\mathcal H^{n-1})$.  Additionally
\begin{equation}
\int_{\partial \Omega} \psi(x, Dv) = \int_{\partial \Omega} \phi(x, Dv +F_{\chi_v})
d\mathcal H^{n-1}
=\ \int_{\partial \Omega} \phi(x, v^+ - v^- +F_{\chi_v})
d\mathcal H^{n-1}.
\label{boundary.meas}\end{equation}
To see this, note that $|Dv|$ can only concentrate on a
set of dimension $n-1$ if that set is a subset of the jump set of $v$,
so \eqref{boundary.meas} follows from standard descriptions of
the jump part of $Dv$.

Now if $u, v\in BV(\R^n)$ satisfy $u=v$ a.e. in $\Omega^c$, then  $\int_{\bar \Omega^c} \varphi(x, Du) = \int_{\bar \Omega^c} \varphi(x, Dv)$. In
addition, $u^+ = v^+$, so
using  \eqref{boundary.meas}
we deduce that
\[
\int_{\R^n}\psi(x, Du) - \int_{\R^n} \psi (x,Dv)
 \ = \ I_\varphi(u; \Omega, u^+) - I_\varphi(v ; \Omega, u^+).
\]
The lemma easily follows from the above equality.
\hfill $\Box$
\\

The next theorem shows super level sets of $\psi$-total variation minimizing functions in $\Omega$ are $\psi$-area minimizing in $\Omega$. 

\begin{theorem} \label{SuperLevelThm}
Let $\Omega \subset \R^n$ be a bounded Lipschitz domain and  $u \in BV(\R^n)$ a $\psi$-total variation minimizing function in $\Omega$.  The super level sets of $u$ are written as
\begin{equation}\label{SuperLevelDef}
    E_\lambda := \left \{ x\in \R^n : u(x) \geq \lambda \right \}.
\end{equation}
Then $E_\lambda$ is $\psi$-area minimizing in $\Omega$.

\end{theorem}

{\bf Proof.}  For a fixed $\lambda \in \R$, let $u_1$ and $u_2$ be as defined in \eqref{u1u2}. Consider $g \in BV(\R^n)$ with $\text{supp}(g) \subset \overline{\Omega}.$ Then
\begin{align*}
    \int_{\Omega} \phi \left( x, D u_1 + F\chi_{\{u\geq \lambda \}} \right) + H u_1 
    &+ \int_{\Omega} \phi \left( x, D u_2 + F\chi_{\{u < \lambda \}} \right) + H u_2 = \int_{\Omega} \phi \left( x,  D u + F \right) + H u\\
    & \leq \int_{\Omega} \phi \left( x,  D (u+g) + F \right) + H (u+g) \\
    &=   \int_{\Omega} \phi \left( x,  D u_1 +D(g\chi_{\{u\geq \lambda\} })+ F\chi_{\{u\geq \lambda \}} \right) + H (u_1+g) \\
    & + \int_{\Omega} \phi \left( x, D u_2 +D(g\chi_{\{u< \lambda\} })+ F\chi_{\{u < \lambda \}} \right) + H u_2 \\
     & \leq   \int_{\Omega} \phi \left( x, D u_1 +D(g\chi_{\{u\geq \lambda\} })+ F\chi_{\{u\geq \lambda \}} \right) + H (u_1+g) \\
    & +\int_{\Omega}\phi \left( x, D(g\chi_{\{u< \lambda\} })\right) + \int_{\Omega} \phi \left( x,  D u_2 + F\chi_{\{u < \lambda \}} \right) + H u_2 \\
    & = \int_{\Omega} \phi \left( x, D (u_1+g) + F\chi_{\{u\geq \lambda \}} \right) + H (u_1+g)\\
    & + \int_{\Omega} \phi \left( x, D u_2 + F\chi_{\{u< \lambda \}} \right) + H u_2.
\end{align*}
This implies
\[ \int_{\Omega} \phi \left( x,  D u_1 + F\chi_{u_1} \right) + H u_1 \leq \int_{\Omega} \phi \left( x,  D (u_1+g) + F\chi_{u_1} \right) + H (u_1+g),\]
for any $g \in BV(\R^n)$ such that $\text{supp}(g) \subset \overline{\Omega}$. By definition, $u_1$ is $\psi$-total variation minimizing.  Using the argument outlined above $\chi_{\epsilon, \lambda}$, as defined in \eqref{chiEpsilon}, is also $\psi$-total variation minimizing.

The boundary of $E_\lambda$ has measure zero for a.e. $\lambda \in \R$, which is represented by
\begin{equation}\label{levelMeasureZero}
     \mathcal{L}^n \left( \{ x \in \Omega : u(x) = \lambda  \} \right) = \mathcal{H}^{n-1} \left( \{ x \in \partial \Omega : u^{\pm}(x) = \lambda  \} \right) =0.  
\end{equation}
 
Thus
$$\chi_{\epsilon, \lambda} \rightarrow \chi_\lambda := \chi_{E_\lambda}  \text{ in } L^1_{\text{loc}}(\R^n), \hspace{0.5cm}  \chi_{\epsilon, \lambda}^{\pm} \rightarrow \chi_{\lambda}^{\pm}  \text{ in } L^1(\partial \Omega ; \mathcal{H}^{n-1}), $$
as $\epsilon \rightarrow 0$.

We apply Lemma \ref{LSC-Lemma} to get
\begin{equation}
P_\psi(\chi_{\lambda}, \Omega) \leq \liminf_{\epsilon \rightarrow 0} P_{\psi}(\chi_{\epsilon,\lambda},\Omega).
\label{wlsc}\end{equation}
It follows from  the $L^1$ convergence of the traces that
\begin{equation}\label{limit1}
I_\varphi(\chi_{\lambda}; \Omega, \chi_{\lambda}^{+}) \le \liminf_{k\to\infty} I_\varphi(\chi_{\epsilon,\lambda}; \Omega, \chi_{\lambda,\epsilon}^{+}).
\end{equation}
For an arbitrary $F\subset \R^n$ with $\chi_{\lambda}=\chi_F$ a.e. in $\Omega^c$, 
\begin{align*}\label{IPhiInequality}
I_\varphi(\chi_{\epsilon,\lambda}; \Omega, \chi_{\epsilon,\lambda}^+) 
&\le
I_\varphi(\chi_F; \Omega, \chi_{\epsilon,\lambda}^+)\\
&\le
I_\varphi(\chi_F; \Omega, \chi_\lambda^+) +\int_{\partial \Omega} \phi(x, \chi_{\lambda}^+-\chi_{\epsilon,\lambda}^{+})
 \ d\mathcal H^{n-1}\\
 &\le
I_\varphi(\chi_F; \Omega, \chi_\lambda^+) +\int_{\partial \Omega} \alpha |\chi_{\lambda}^+-\chi_{\epsilon,\lambda}^{+}|
 \ d\mathcal H^{n-1}\\
&\le
I_\varphi(\chi_F; \Omega, {\chi_{\lambda}}^+) +C\int_{\partial \Omega}
|\chi_{\lambda}^+-\chi_{\epsilon,\lambda}^{+}|\ d\mathcal H^{n-1}.
\end{align*}
The inequality that follows is justified by the above, \eqref{limit1},
and $\chi_{\epsilon,\lambda}^+ \rightarrow \chi_{\lambda}^+$ in $L^1(\partial \Omega;\mathcal H^{n-1})$,

\[I_\varphi(\chi_\lambda;\Omega, \chi_{\lambda}^+) \le 
I_\varphi(\chi_F;\Omega, \chi_\lambda^+).\]
This establishes that $E_{\lambda}$ is $\phi$-area minimizing in $\Omega$.

If $\lambda$ does not satisfy (\ref{levelMeasureZero}), then there exists an increasing sequence $ \lambda_k  $ that converges to $\lambda$ and satisfies (\ref{levelMeasureZero}) for each $k$.  In which case,
$$\chi_{\lambda_k} \rightarrow \chi_\lambda \text{ in } L^1_{\text{loc}}(\R^n), \hspace{0.5cm}  \chi_{\lambda_k}^{\pm} \rightarrow \chi_{\lambda}^{\pm}  \text{ in } L^1(\partial \Omega ; \mathcal{H}^{n-1}). $$
Thus, by  Lemma \ref{Relaxed-Lemma}, $E_\lambda$ is $\psi$-area minimizing in $\Omega$.  \hfill  $\Box$
\\

It remains to lay out a few more definitions which would play a key role in the proof of our main result in this section. Let
$$ 
BV_f(\Omega) := \left\{ u \in BV(\Omega): \lim_{r \rightarrow 0} \esssup_{y \in \Omega, |x-y|<r} |u(y)-f(y)|=0 \text{ for } x \in \partial \Omega \right\}. 
$$
For any measurable set $E$, consider

\[E^{(1)}:= \{x\in \R^n :
\lim_{r\to 0} \frac {\calH^n(B(r,x)\cap E)}{\calH^n(B(r))} = 1 \}.\]

\begin{definition} \label{BarrierCondition}
Let $\Omega \subset \R^n$ be a bounded Lipschitz domain. We say that  $\Omega$ satisfied the barrier condition if for every $x_0 \in \partial \Omega$ and $\epsilon >0$ sufficiently small,  $V$ minimizes $P_\psi ( \cdot ; \R ^n)$ in
\begin{equation}\label{BD-cond}
\{ W \subset  \Omega: W \setminus B(\epsilon, x_0 ) = \Omega \setminus B(\epsilon, x_0 )    \},
\end{equation}
implies
$$
\partial V^{(1)} \cap \partial \Omega \cap B(\epsilon, x_0)= \emptyset.
$$
\end{definition}

Intuitively speaking, (\ref{BD-cond}) means that at any point $x_0\in \partial \Omega$ one can decrease the $\psi$-perimeter of $\Omega$ by pushing the boundary inwards,  

\begin{lemma} \label{MainLemma}
    Suppose $\Omega \subset \R^n$ is a bounded Lipschitz domain satisfying the barrier condition, and $E \subset \R^n$ minimizes $P_\psi(\cdot ;\Omega)$.  Then
    $$
    \left\{ x \in \partial \Omega \cap \partial E^{(1)} : B(\epsilon, x) \cap \partial E^{(1)} \subset \overline{\Omega} \text{ for some } \epsilon > 0 \right\} = \emptyset.
    $$
\end{lemma}
{\bf Proof.} We proceed by contradiction. Suppose there exists $x_0\in \partial \Omega \cap \partial E^{(1)}$ such that  $B(\epsilon, x_0) \cap \partial E^{(1)}\subset \bar{\Omega}$ for some $\epsilon>0$. Then $\tilde{V}= E\cap \Omega$ is a minimizer of $P_{\psi}(\,\cdot\,; \R^n)$ in (\ref{BD-cond}), and 
\[
x_0\in \partial {\tilde{V} }^{(1)}\cap \partial \Omega \cap B(\epsilon, x_0)\neq  \emptyset. 
\]
This is inconsistent with the barrier condition (\ref{BD-cond}). 
\hfill $\Box$

\medskip

Finally, we are ready to prove the main existence results of the this section.

\begin{theorem} \label{MainThm}
Consider $\psi:\R ^n \times \R ^n \rightarrow \R$ as defined in \eqref{AlteredIntegrand} and a bounded Lipschitz domain $\Omega \subset \R ^n$. Let $||H||_{L^{\infty}(\overline{\Omega})}$ be small enough that Proposition \ref{PropLast} holds. If $\Omega$ satisfies the barrier condition with respect to $\psi$, then for every $f \in C(\partial \Omega)$ the minimization problem \eqref{functionalMain} has a minimizer $u$ in $BV(\Omega)$ with $u|_{\partial \Omega} \leq f $.
\end{theorem}

{\bf Proof.}  For a given $f \in C(\partial \Omega)$, it can be extended to $f \in C(\Omega ^c)$.  Furthermore, we can assume $f \in BV(\R^n)$ since every $\calH^{n-1}$
integrable function on $\Omega$ is the trace of some
(continuous) function in $BV(\Omega^c)$.  Let
\[
\calA_f := \{v\in BV(\R^n): \ \ v=f \ \ \hbox{on}\ \ \Omega^c \},
\]
where any element $v$ of $BV_f(\Omega)$
is the restriction to $\Omega$ of a unique element of 
$\calA_f$. Then $\int_{\R^n} \psi(x,v)$ has as a minimizer $u\in \calA_f$, in view of Proposition \ref{PropLast}.  

Next we prove that $u|_{\partial \Omega}\leq f$. Suppose this is not the case, then there is an $x\in \partial \Omega$
and $\delta>0$ such that
\begin{equation} 
  \esssup_{y\in \Omega, |x-y|<r}\big(u(y) - f(x)) \ge \delta
\label{ess.alt}\end{equation}
for every $r>0$.
First, suppose that the latter condition holds.  For $E := E_{ f(x) + \delta/2}$ we have that $x\in  \partial E^{(1)}$, justified by the second alternative of \eqref{ess.alt} and the continuity of $f$.  Note that Theorem \ref{SuperLevelThm} implies $E$ is $\psi$-area minimizing in $\Omega$.  This there exists  $\epsilon >0$ such that $u< f(x) + \delta/2$ in $B(\e, x)\setminus \Omega$, since $u\in \calA_f$ and $f$ is continuous in $\Omega^c$.  However, Lemma \ref{MainLemma} shows that this is impossible.   \hfill $\Box$   \\ \\


\begin{thebibliography}{99}


\bibitem {Al} {\sc G. Alberti}, {\em A Lusin type theorem for gradients}, J. Funct. Anal., Vol. \textbf{100} (1991), pp. 110-118.


\bibitem{AB} {\sc M. Amar, G. Bellettini}, {\em A notion of total variation depending on a metric with discontinuous coefficients}, Annales de l'institut Henri Poincar\'{e}(C) Analyse non linéaire \textbf{11} (1994),  91-133.



\bibitem{An} {\sc G. Anzellotti}, {\em Pairings between measures and bounded functions and compensated compactness}, 
Ann. Mat. Pura Appl. (4) 135 (1983), 293-318 (1984). 

\bibitem{And4} {\sc F. Andreu-Vaillo, V.  Caselles, J. M. Maz\'{o}n}, 
 {\em Parabolic quasilinear equations minimizing linear growth functionals}, Progress in Mathematics, 223. Birkhäuser Verlag, Basel, 2004.


\bibitem{B} {\sc Z.M. Balogh}, \emph{Size of characteristic sets and functions with prescribed gradient}. J. Reine Angew.
Math. 564 (2003), 63-83.

\bibitem{Bousq} {\sc P. Bousquet}, {\em Boundary continuity of solutions to a basic problem in the calculus of variations}, Adv. Calc. Var. \textbf{3} (2010), 1-27. 

\bibitem{BC} {\sc P. Bousquet, F. Clarke}, {\em Local Lipschitz continuity of solutions to a problem in the calculus of variations}, J. Differential Equations \textbf{243} (2007), 489–503. 

\bibitem{C} {\sc A. Cellina}, {\em On the bounded slope condition and the validity of the Euler Lagrange equation}, SIAM J. Control Optim. \textbf{40} (2001/02), 1270–1279 (electronic). 


\bibitem{CHP} {\sc J.-H. Cheng,  J.-F. Hwang}, \emph{Properly embedded and immersed minimal surfaces in the Heisenberg group}. Bull. Aus. Math. Soc. 70 (2004), 507-520.

\bibitem{CH} {\sc J.-H. Cheng,  J.-F. Hwang}, \emph{Uniqueness of generalized p-area minimizers and integrability of a horizontal normal in the Heisenberg group}. 
Calc. Var. Partial Differential Equations 50 (2014), no. 3-4, 579-597.

\bibitem{CHMY} {\sc J.-H. Cheng,  J.-F. Hwang, A. Malchiodi, P. Yang},  \emph{Minimal surfaces in pseudohermitian geometry}. Annali della Scuola Normale Superiore di Pisa, Classe di Scienze 4(5) (2005), 129-177.

\bibitem{CHY} {\sc J.-H. Cheng,  J.-F. Hwang, A. Malchiodi, P. Yang},
\emph{Existence and uniqueness for p-area minimizers in the Heisenberg group}.  Math. Ann. 337 (2007), no. 2, 253-293.

\bibitem{Clar} {\sc F. Clarke}, {\em Continuity of solutions to a basic problem in the calculus of variations}, Ann. Sc. Norm. Super. Pisa
Cl. Sci. (5), \textbf{4} (2005), 511–530. 

\bibitem{DLPT} {\sc S. Don,  L. Lussardi, A. Pinamonti, G. Treu},
\emph{Lipschitz minimizers for a class of integral functionals under the bounded slope condition}.  Nonlinear Analysis, Theory, Methods and Applications, 216 (2022), 112689.

\bibitem{ET}  {\sc I. Ekeland, R. T\'{e}mam}, \emph{Convex analysis and variational problems}, North-Holland-Elsevier, 1976. 

\bibitem{FT} {\sc A. Fiaschi, G. Treu}, {\em The bounded slope condition for functionals depending on $x, u,$ and $\nabla u$}, SIAM J. Control
Optim., \textbf{50} (2012), 991–1011. 


\bibitem{FSS} {\sc B. Franchi, R. Serapioni, F. Serra Cassano}, \emph{Rectifiability and perimeter in the Heisenberg
group}. Math. Ann. 321, 479-531 (2001).
\bibitem{GN} N. Garofalo,  D.-M Nhie, \emph{ Isoperimetric and Sobolev inequalities for Carnot-Caratheodory
spaces and the existence of minimal surfaces}. Comm. Pure Appl. Math. 49, 1081-1144 (1996).


\bibitem{G} {\sc E. Giusti}, {\em Minimal Surfaces and Functions of Bounded
Variations}, Birkh\"{a}user, Boston, 1984.



\bibitem{Gorny} {\sc W. G\'{o}rny}, {\em Planar least gradient problem: existence, regularity and anisotropic case}, https://arxiv.org/abs/1608.02617. 


\bibitem{HMN} {\sc N. Hoell, A. Moradifam, A. Nachman}, {\em Current Density Impedance Imaging with an Anisotropic Conductivity in a Known Conformal Class}, SIAM J. Math. Anal., 46 (2014), 3969-3990.

\bibitem{JMN} {\sc R.L. Jerrard, A. Moradifam, A. Nachman}, {\em Existence and uniqueness of minimizers of general least gradient problems}, J. Rein Angew. Math., 734 (2018), 71-97.


\bibitem{LM} {\sc L. Lussardi, E.Mascolo}, {\em A uniqueness result for a class of non strictly convex variational problems}, J. Math. Anal. Appl. \textbf{446} (2017), no. 2, 1687–1694. 

\bibitem{MT} {\sc C. Mariconda, G. Treu}, {\em Existence and Lipschitz regularity for minima}, Proc. Amer. Math. Soc. \textbf{130} (2002), 395–404 (electronic). 

\bibitem{MT2} {\sc C. Mariconda, G. Treu}, {\em Lipschitz regularity for minima without strict convexity of the Lagrangian}, J. Differential Equations \textbf{243} (2007), 388–413. 

\bibitem{MT3} {\sc C.  Mariconda, G. Treu}, {\em Local Lipschitz regularity of minima for a scalar problem of the calculus of variations}, Commun. Contemp. Math. \textbf{10} (2008), 1129–1149. 

\bibitem{Mazon2} {\sc J. M. Maz\'{o}n}, {\em The Euler–Lagrange equation for the Anisotropic least gradient problem}, Nonlinear Analysis: Real World Applications \textbf{31} (2016) 452-472. 

\bibitem{MRD} {\sc J. M. Maz\'{o}n, J.D. Rossi, S.S. De Le\'{o}n }, {\em Functions of Least Gradient and 1-Harmonic Functions}, Indiana University Mathematics Journal \textbf{63} (2013) (4): 1067-1084.

\bibitem{Mo} {\sc A. Moradifam}, {\em Existence and structure of minimizers of  least gradient problems},  Indiana University Mathematics Journal \textbf{63} (2014), no. 6, 1819-1837. 

\bibitem{Mo1} {\sc A. Moradifam}, {\em Least gradient problems with Neumann boundary condition}, J. Differential Equations 263 (2017), no. \textbf{11}, 7900-7918.

\bibitem{MNT}
{\sc A. Moradifam, A. Nachman, and A. Timonov}, {\em A convergent algorithm for the hybrid problem of reconstructing conductivity from minimal interior data}, Inverse Problems,  {\bf 28} (2012) 084003.

\bibitem{MNTa_SIAM} {\sc A. Moradifam, A. Nachman, and A. Tamasan}, {\em Conductivity imaging from one interior measurement in the presence of perfectly conducting and insulating inclusions}, SIAM J. Math. Anal., {\bf 44} (2012)  (6), 3969-3990.

\bibitem{MR} {\sc A. Moradifam, A. Rowell}, {\em Existence and structure of P-area minimizing surfaces in the Heisenberg group}, Journal of Differential Equations, {\bf 342} (2023), 325-342. 

\bibitem{NTT07}
{\sc A. Nachman, A. Tamasan, and A. Timonov}, {\em Conductivity
imaging with a single measurement of boundary and interior data},
Inverse Problems, {\bf 23} (2007), pp. 2551--2563.

\bibitem{NTT08}
{\sc A. Nachman, A. Tamasan, and A. Timonov}, {\em Recovering the
conductivity from a single measurement of interior data}, Inverse
Problems, {\bf 25} (2009) 035014 (16pp).


\bibitem{NTT10}{\sc A. Nachman, A. Tamasan, and A. Timonov}, {\em Reconstruction of Planar Conductivities in Subdomains from Incomplete Data},
SIAM J. Appl. Math. {\bf 70}(2010), Issue 8, pp. 3342--3362.

\bibitem{NTT11}{\sc A. Nachman, A. Tamasan, and A. Timonov}, {\em Current density impedance imaging}, Tomography and inverse transport theory, 135-149, Contemp. Math. {\bf 559}, AMS, 2011.

\bibitem{P} {\sc S.D. Pauls},\emph{ Minimal surfaces in the Heisenberg group}. Geometric Dedicata, 104 (2004), 201-231.

\bibitem{PSTV} {\sc A. Pinamonti, F.  Serra Cassano, G. Treu, D. Vittone}, \emph{BV minimizers of the area functional in the Heisenberg group under the bounded slope condition}. Ann. Sc. Norm. Super. Pisa Cl. Sci. (5) 14 (2015), no. 3, 907-935.

\bibitem{sternberg_ziemer92} {\sc P. Sternberg, G. Williams, and W. P. Ziemer}, {\em Existence, uniqueness and regularity for functions of least gradient}, J. Rein Angew. Math. {\bf 430} (1992), 35-60.

\bibitem{sternberg_ziemer}{\sc P. Sternberg and W. P. Ziemer}, {\em Generalized motion by
curvature with a Dirichlet condition}, J. Differ. Eq., {\bf 114}(1994), pp. 580--600.

\bibitem{sternbergZiemer93} {\sc P. Sternberg and W. P. Ziemer}, {\em The Dirichlet problem for
functions of least gradient. Degenerate diffusions} (Minneapolis, MN,
1991), 197--214, in IMA Vol. Math. Appl., \textbf{47}, Springer, New York, 1993.

\bibitem{ST} {\sc G. S. Spradlin and A. Tamasan}, {\em Not all traces on the circle come from functions of least gradient in the disk}, Indiana University Mathematics Journal \textbf{63} (2014), no. 6, 1819-1837. 












\end{thebibliography}
\end{document}